\documentclass[10pt]{amsart}

\usepackage{amssymb}
\usepackage[all]{xy}
\usepackage{mathrsfs}
\usepackage{enumerate}

\usepackage{bbm}

\newcommand{\rar}[1]{\stackrel{#1}{\longrightarrow}}

\newcommand{\bC}{{\mathbb C}}
\newcommand{\bD}{{\mathbb D}}

\newcommand{\bZ}{{\mathbb Z}}

\newcommand{\cA}{{\mathcal A}}
\newcommand{\cB}{{\mathcal B}}
\newcommand{\cC}{{\mathcal C}}

\newcommand{\cF}{{\mathcal F}}
\newcommand{\cG}{{\mathcal G}}

\newcommand{\cI}{{\mathcal I}}

\newcommand{\cR}{{\mathcal R}}

\newcommand{\cT}{{\mathcal T}}

\newcommand{\nc}{\newcommand}

\nc\wh{\widehat}

\nc\on{\operatorname}

\nc\Gr{\on{Gr}}

\nc\Fl{\on{Fl}}

\newcommand{\limto}{{\displaystyle\lim_{\longrightarrow}}}
\newcommand{\rightlim}{\mathop{\limto}}

\newcommand{\leftlim}{\mathop{\displaystyle\lim_{\longleftarrow}}}
\newcommand{\limfromn}{\leftlim\limits_{\raise3pt\hbox{$n$}}}
\newcommand{\limton}{\rightlim\limits_{\raise3pt\hbox{$n$}}}

\newcommand{\rightlimit}[1]{\mathop{\lim\limits_{\longrightarrow}}\limits%
                    _{\raise3pt\hbox{$\scriptstyle #1$}}}

\newcommand{\leftlimit}[1]{\mathop{\lim\limits_{\longleftarrow}}\limits%
                    _{\raise3pt\hbox{$\scriptstyle #1$}}}

\newcommand{\iso}{\buildrel{\sim}\over{\longrightarrow}}
\newcommand{\mono}{\hookrightarrow}

\DeclareMathOperator{\Res}{{Res}}

\newtheorem{Th}{Theorem}
\newtheorem{pr}{Proposition}[section]
\newtheorem{lm}[pr]{Lemma}

\newtheorem{cor}[pr]{Corollary}

\theoremstyle{definition}

\newtheorem{rem}[pr]{Remark}

\numberwithin{equation}{section}

\begin{document}

\title{On the derived DG functors}

\author{Vadim Vologodsky}




\begin{abstract}
Assume that abelian categories ${\cA}$,  ${\cB}$ over a field admit countable direct limits and that these limits are exact.  Let $\cF:   D^+_{dg}({\cA}) \to D^+_{dg}( {\cB})$ be a DG quasi-functor  such that  the functor $Ho(\cF):  D^+({\cA}) \to D^+({\cB})$ carries  $D^{\geq 0}({\cA})$ to $D^{\geq 0}({\cB})$
 and such that, for every    $i>0$,  the functor $H^i \cF: \cA \to \cB$ is
  effaceable.  We prove that $\cF$  is canonically isomorphic to the right derived DG functor $RH^0(\cF)$. 
We also prove a similar result for bounded derived DG categories and a formula that expresses Hochschild cohomology of  the categories $ D^b_{dg}({\cA})$, $ D^+_{dg}({\cA}) $ as the $Ext$ groups in the abelian category of left exact functors $\cA \to Ind \cA$ .  The proofs are based
on a description of Drinfeld's category of quasi-functors as the derived category of a category of sheaves. 
\end{abstract}

\maketitle


\section{Main results}\label{mr}



Let $\cA$ and $\cB$ be abelian categories, and let
$$RF_{tri}: D^{+}({\cA}) \to D^{+}({\cB})$$
be the right derived functor of some left exact functor $F: \cA \to \cB$. Then, the corresponding  cohomological $\delta$-functor $R^*F= H^* RF_{tri}:  \cA \to \cB$  has the following property:  the functor $H^iRF_{tri}$ is $0$ for $i<0$,  effaceable for $i>0$, and
$H^0RF_{tri}$  is isomorphic to $F$. 
Conversely, according to a result of Grothendieck (\cite{g}), every  cohomological $\delta$-functor $T^*: \cA \to \cB$  satisfying the above property is canonically isomorphic to the right derived functor $R^*F$. 
  The purpose of this paper is to extend this
extremely useful characterization of $R^*F$ to the derived category level.  
Unfortunately, Verdier's notion of triangulated functor seems  too poor to allow such a simple characterization of the derived functors.  In order to get a meaningful statement one has to consider
triangulated functors with some kind of enrichment.  Arguably the most useful notion here is the one of {\it DG quasi-functor} (or essentially equivalent  notion of $A_{\infty}$-functor).   Indeed, works of Keller and Drinfeld (\cite{k2}, \cite{dri}) provide a canonical
DG enhancement  $ D^{+}_{dg}({\cA})$ of  Verdier's triangulated derived category.  Roughly, a DG  quasi-functor $\cF:  D^{b}_{dg}({\cA}) \to D^{b}_{dg}({\cB}) $ is a diagram of the form  
\begin{equation}\label{qdgf}
       D^{+}_{dg}({\cA}) \stackrel{S}{\longleftarrow}\cC \rar{G} D^{+}_{dg}({\cB}).
       \end{equation}
       where $\cC$ is a DG category, $S$ and $G$ are DG functors, and, in addition,  $S$ is a homotopy equivalence.
  Every quasi-functor  (\ref{qdgf}) yields
a triangulated functor $Ho(\cF): D^{+}({\cA}) \to D^{+}({\cB})$, but the converse is not true in general. Nevertheless,  many of the natural triangulated functors come together with a DG enhancement. For example, the triangulated derived functor $RF$ can be canonically 
promoted to a DG quasi-functor (\cite{dri} \S 5). The main result of this paper states that under certain mild assumptions on abelian categories $\cA$ and $\cB$ the DG quasi-functors isomorphic to the DG derived ones are precisely the DG quasi-functors satisfying
Grothendieck's condition above.  
 To state the result we need to introduce a bit of notation.

Let $k$ be a commutative ring.  Denote by $Mod(k)$ the category of $k$-modules. We shall say that a $k$-linear category \footnote {{\it i.e.}, a category enriched over $Mod(k)$.} is $k$-flat if, for every two objects $X,Y$,  the $k$-module $Hom(X,Y)$ is flat.  
Given a $k$-linear exact category  ${\cA}$  we denote by   $D^{b}_{dg}({\cA})$ the corresponding bounded derived DG category over $k$.
This is the DG quotient (\cite{dri}) of the DG category $C^b_{dg}({\cA})$ of bounded complexes by the subcategory of acyclic ones (\cite{n}, \S 1). The homotopy category of $D^{b}_{dg}({\cA})$ is the triangulated derived category $D^{b}({\cA})$ as defined in (\cite{n}).  
Let ${\cB}$ be another $k$-linear abelian category, $ D^{b}_{dg}({\cB})$ the corresponding bounded derived DG category, 
and let  $\cT(D^{b}_{dg}({\cA}), D^{b}_{dg}({\cB}) )   $ be the triangulated category of  DG quasi-functors  $\cF:  D^{b}_{dg}({\cA}) \to D^{b}_{dg}({\cB}) $   (\cite{dri}, \S 16.1).  Given such $\cF$  and an integer $i$  we denote by $H^i\cF: \cA \to \cB$ the composition
$$ \cA \to D^{b}_{dg}({\cA}) \stackrel{\cF}{\longrightarrow} D^{b}_{dg}({\cB}) \stackrel{H^i}{\longrightarrow} \cB. $$

 \begin{Th}\label{bounded}
  Let $\cA$ be a small $k$-flat exact idempotent complete category \footnote{An additive category is called  idempotent complete if any its morphism $p: X\to X$ such that $p\circ p=p $ is the projection on a direct summand of a decomposition $X\simeq Y\oplus Z$.}
   and $\cB$ a  small abelian $k$-linear category. 
  \begin{enumerate}
\item{Assume that a DG quasi-functor  $$\cF:  D^{b}_{dg}({\cA}) \to D^{b}_{dg}({\cB}) $$ 
 has the following property:

(P) The functor $H^i \cF: \cA \to \cB$ is $0$ for every $i<0$ and
  effaceable  \footnote{That is,  for every object $X\in \cA$, there exists an admissible monomorphism  $X\hookrightarrow Y$ such that 
the induced morphism   $H^i \cF(X) \to  H^i \cF(Y)$ is $0$.} for every $i>0$.

Then the functor $F:= H^0 \cF : \cA \to \cB$ is left exact, has a right derived DG quasi-functor (\cite{dri} \S 5)  
         $$RF:  D^{b}_{dg}({\cA}) \to D^{b}_{dg}({\cB}), $$ 
         and  there is a unique isomorphism $\cF \simeq RF$ such that the induced automorphism  $F= H^0(\cF) \simeq H^0(RF)= F$ equals $Id$. Conversely,  the right derived DG quasi-functor of any  left exact functor $F:  \cA \to \cB$ satisfies property (P). } 
         
\item{For every two DG quasi-functors $\cF,\cG  \in \cT(D^{b}_{dg}({\cA}), D^{b}_{dg}({\cB}) ) $ satisfying property (P)
and every $i<0$, we have
$$Hom_{\cT(D^{b}_{dg}({\cA}), D^{b}_{dg}({\cB}) )}(\cF, \cG[i])=0,$$
$$Hom _{\cT(D^{b}_{dg}({\cA}), D^{b}_{dg}({\cB}) )}(\cF, \cG) = Hom_{Fct(\cA, \cB)}(H^0\cF, H^0\cG).$$
Here $Fct(\cA, \cB)$ denotes the category of all $k$-linear functors $\cA \to \cB$.  }

\end{enumerate}

 \end{Th}
  \begin{rem}\label{tri} I do not know if the analogous statement holds for merely triangulated functors. 
 \end{rem}
  \begin{rem} It is likely that  the $k$-flatness assumption on $\cA$ is unnecessary. However, I can not prove this. 
 \end{rem}

 We have a similar result for bounded from below derived DG categories. If $\cA$ is a $k$-linear abelian category we will write $D^{+}_{dg}({\cA})$ for the bounded from below derived DG category of $\cA$ and  $D^{+}({\cA})$ for the corresponding triangulated category.
 Let $ D^{\geq n}({\cA}) $ be the full subcategory of $D^{+}({\cA})$ that consists of  complexes with trivial cohomology in degrees less then $n$.  We say that  a DG quasi-functor $$\cF:  D^{+}_{dg}({\cA}) \to D^{+}_{dg}({\cB}) $$ has property ($P'$) if
 
 ($P'$) The functor $Ho(\cF)$ takes every object of  the category $D^{\geq 0}({\cA}))$ to an object of  $D^{\geq 0}({\cB})$ and, for every    $i>0$,  the functor $H^i \cF: \cA \to \cB$ is
  effaceable.  
 
 \begin{Th}\label{unbounded} Let $k$ be a field and let $\cA$, $\cB$ be small abelian $k$-linear categories.   Assume that both categories are closed under countable direct limits and that these limits are exact. 
   \begin{enumerate}
   \item{Let  $\cF  \in  \cT(D^{+}_{dg}({\cA}), D^{+}_{dg}({\cB}) ) $ be  a DG quasi-functor satisfying property ($P'$) and $F:= H^0 \cF: \cA \to \cB$.  
   The functor $F$ admits a right derived DG quasi-functor   
         $RF:  D^{+}_{dg}({\cA}) \to D^{+}_{dg}({\cB})$ and there is a unique isomorphism $\cF \simeq RF$ such that the induced automorphism  $F= H^0(\cF) \simeq H^0(RF)= F$ equals $Id$. 
            Conversely,  a right derived DG quasi-functor of any  left exact functor $F:  \cA \to \cB$ satisfies property ($P'$). } 
      \item{For every two DG quasi-functors $\cF,\cG  \in  \cT(D^{+}_{dg}({\cA}), D^{+}_{dg}({\cB}) ) $ satisfying property ($P'$)
and every $i<0$, we have
$$Hom_{ \cT(D^{+}_{dg}({\cA}), D^{+}_{dg}({\cB}) ) }(\cF, \cG[i])=0,$$
$$Hom _{ \cT(D^{+}_{dg}({\cA}), D^{+}_{dg}({\cB}) ) }(\cF, \cG) = Hom_{Fct(\cA, \cB)}(H^0\cF, H^0\cG).$$  }
         \end{enumerate}

 \end{Th}
 The main ingredient of the proof of Theorem \ref{unbounded} is the following construction. Let  $Sh(\cA^{o} \otimes _k \cB) $ be the category of $k$-linear contravariant functors $\cA^{o} \otimes _k \cB \to Mod(k)$ that are left exact with respect to both arguments.
 Every $k$-linear left exact functor $F: \cA \to \cB$ yields $s(F)\in  Sh(\cA^{o} \otimes _k \cB) $:
      $$s(F)(X\otimes X')= Hom_{\cB}(X', F(X)).$$
Let $\cT^+ \subset \cT(D^{+}_{dg}({\cA}), D^{+}_{dg}({\cB}))$ be the full triangulated subcategory whose objects are quasi-functors $\cF$ such that
$ Ho(\cF)( D^{\geq 0}({\cA})) \subset  D^{\geq n}({\cB}) $ for some $n$. Using key Lemma \ref{key}  we construct a fully faithful embedding 
\begin{equation}\label{idea}
\cT^+ \mono D(Sh(\cA^{o} \otimes _k \cB) )
\end{equation}
that carries every DG quasi-functor $\cF$ satisfying property ($P'$) to  $s(F)\in Sh(\cA^{o} \otimes _k \cB) \subset   D(Sh(\cA^{o} \otimes _k \cB) ).$ 
\begin{rem}   In (\cite{t},  Th. 8.9), To\"en gave an analogous description of the category of quasi-functors between the derived DG categories of (quasi)-coherent sheaves. 
   \end{rem}
As another application of (\ref{idea}) we compute the Hochschild cohomology of a derived DG category. Recall (see, e.g. \cite{k1}, \S 5.4,  \cite{t}, \S 8.1) that  the Hochschild cohomology  of a DG category $\cC$ can be interpreted as
\begin{equation}\label{hoch}
HH^i(\cC, \cC)= Hom_{\cT(\cC, \cC)}(Id_{\cC}, Id_{\cC}[i]).
\end{equation}
The composition in $\cC$ makes $HH^*(\cC, \cC)$ a graded commutative algebra over $k$.
 \begin{Th}\label{Hochschild}
Let $k$ be a field, and let $\cA$ be a small abelian $k$-linear category.  There is an isomorphism of algebras 
\begin{equation}\label{hoch1}
HH^*(D^{b}_{dg}({\cA}), D^{b}_{dg}({\cA}) )\simeq Ext^*_{Sh(\cA^{o} \otimes _k \cA)}(s(Id_{\cA}), s(Id_{\cA})).
\end{equation}
 If, in addition,  $\cA$ is closed under countable direct limits and that these limits are exact, we have
 \begin{equation}\label{hoch2}
 HH^*(D^{+}_{dg}({\cA}), D^{+}_{dg}({\cA}) )\simeq Ext^*_{Sh(\cA^{o} \otimes _k \cA)}(s(Id_{\cA}), s(Id_{\cA})).
 \end{equation}
 \end{Th}
 \begin{rem}  This is a remarkable phenomenon the Hochschild cohomology does not change we ``enlarge'' the DG category.  A similar result, that  the Hochschild cohomology of a small DG category coincides with 
  the Hochschild cohomology of its DG ind-completion, is due to To\"en  (\cite{t}, \S 8).  An analogous statement for Grothendieck abelian categories was proved by Lowen and  Van den Bergh (\cite{lv}).
  \end{rem}
 
  \begin{rem}  The category $Sh(\cA^{o} \otimes _k \cA)$ has a tensor structure that extends the tensor structure on the category of left exact endofunctors $\cA \to \cA$ given by the composition.  
  This  can be used to promote (\ref{hoch1}), (\ref{hoch2}) to isomorphisms of {\it Gerstenhaber algebras}  (see, e.g. \cite{k1}, \S 5.4). 
   \end{rem}

 {\bf Notation.}   Given a category $\cC$ we denote by $\cC^o$ the opposite category. If $\cC$ is a DG category  we will write $Ho\, \cC$ for the corresponding homotopy category  (\cite{dri}, \S 2.7). For example,  $Ho\,C(Mod(k))$ denotes the homotopy  category of complexes of $k$-modules.  The derived category of right DG modules over a DG category 
 $\cC$ will be denoted by $\bD(\cC)$  (\cite{dri}, \S 2.3) \footnote{Drinfeld's notation for this category is $D(\cC)$. We use a different notation to avoid a possible confusion with Verdier's  derived category of an abelian category $\cC$ that is denoted by $D(\cC)$.}.
 We will write  $\underrightarrow{\cC}$ for the DG category of semi-free right DG modules over $\cC$ (\cite{bv}, 1.6.1). We have a canonical equivalence of triangulated categories  $Ho\underrightarrow{\cC}\iso \bD(\cC)$ (\cite{bv}, 1.6.4). For DG categories 
 $\cC$, $\cC'$ we denote by $\cT(\cC, \cC')$ the category of DG quasi-functors  (\cite{dri}, \S 16.1). If $\cC'$ is a pretriangulated  (\cite{dri}, \S 2.4) $\cT(\cC, \cC')$  has a canonical structure of triangulated category.  If $\cF \in \cT(\cC, \cC')$ we will
 write $Ho(\cF)$ for the corresponding functor between the homotopy categories.  The expression ``direct limit'' always means ``filtrant direct limit'' (\cite{ks}, \S 3).

{\bf Acknowledgements.} I would like to thank Sasha Beilinson, Bernhard Keller, and Dima Orlov for helpful conversations related to the subject of this paper.  My deep thanks due to the referee for his numerous remarks and for pointing out an error in a preliminary 
version of this paper. This research was partially supported by NSF grant  DMS-0901707.  Writing this paper I had in mind an application to the theory of Voevodsky's motives (\cite{v}). However, I believe that the main result explained here is interesting on its own ground. 

\section{Proofs}\label{proofs}
{\bf Proof of theorem \ref{bounded}.} Let $\cT^+ \subset \cT:= \cT(D^{b}_{dg}({\cA}), D^{b}_{dg}({\cB}) )$ be the full triangulated subcategory whose objects are quasi-functors $\cF$ such that $H^i \cF =0 $ for sufficiently small $i$. 
To prove the Theorem,  we shall construct (in Lemma \ref{key} below) a fully faithful  embedding of  $\cT^+$  into the derived category of a certain abelian category $Sh(\cA^{o} \otimes _k \cB) $ that takes every functor
$\cF \in \cT^+$ satisfying property (P)  to an object of the heart   $Sh(\cA^{o} \otimes _k \cB) \subset D(Sh(\cA^{o} \otimes _k \cB))$.

Under our flatness assumption on $\cA$, the category  ${\cT}$ is a full subcategory of the derived category  $\bD( D^{b}_{dg}({\cA})^{o}\otimes_k D^{b}_{dg}({\cB}))$ of  right DG modules over 
$D^{b}_{dg}({\cA})^{o}\otimes_k D^{b}_{dg}({\cB})$ that  consists of all $M \in \bD( D^{b}_{dg}({\cA})^{o}\otimes_k D^{b}_{dg}({\cB}))$ such that, for every $X$ in $D^{b}_{dg}({\cA})^{o}$, the module $M(X) \in  \bD( D^{b}_{dg}({\cB}))$ belongs to the essential image of the Yoneda embedding
 $ D^{+}_{dg}({\cB}) \to \bD( D^{b}_{dg}({\cB}))$   (\cite{dri}, \S 16.1).   

Consider the restriction functor
 $$ \bD( D^{b}_{dg}({\cA})^{o}\otimes_k D^{b}_{dg}({\cB})) \rar{\beta} \bD( \cA ^{o}\otimes_k  \cB )$$
  induced by the DG quasi-functor $  \cA ^{o}\otimes_k  \cB  \to D^{b}_{dg}({\cA})^{o}\otimes_k D^{b}_{dg}({\cB})$.
 By definition, the triangulated category   $\bD( \cA ^{o}\otimes_k  \cB )$  is the derived category of the abelian category $ PSh: =PSh( \cA ^{o} \otimes_k  \cB )$  of $k$-linear presheaves {\it i.e.}, the category of $k$-linear contravariant functors $\cA^{o} \otimes _k \cB \to Mod(k)$.  
 Consider a Grothendieck topology on $ \cA^{o} \otimes _k \cB $ whose covers are maps of the form $f\otimes g : Y\otimes Y' \to X \otimes X'$, where $X, Y \in \cA^{o}$, $X', Y' \in \cB$, and $f: Y\to X$, $g: Y' \to X'$  are  admissible
epimorphisms \footnote{By definition,  admissible epimorphisms $Y\to X$ in $\cA^{o}$ are admissible monomorphisms $X\to Y$ in $\cA$.} {\it i.e.,} a sieve $\cC$ over $X\otimes X'$ is a covering sieve if there exist  $f: Y\to X$, $g: Y' \to X'$ as above such that
$ Y\otimes Y' \rar{f\otimes g}  X \otimes X' \in \cC$.
  The axioms of Grothendieck topology (see, e.g. \cite{ks}, \S 16.1)  are immediate except for the one which is the following statement: for every cover  $ Y\otimes Y' \rar{f\otimes g}  X \otimes X'$ and every morphism $ Z\otimes Z' \rar{\phi}  X \otimes X'$  there exists
  a cover $ T\otimes T' \rar{p\otimes q}  Z \otimes Z'$  and a morphism  $ T\otimes T' \rar{\psi}  Y \otimes Y'$  such that $(f\otimes g) \circ \psi=   \phi \circ (p\otimes q)$, which is a consequence of  the base change axiom of exact category (\cite{q}, \S 2). 
 Let $Sh:= Sh( \cA^{o} \otimes _k \cB )$ be the subcategory of  $ PSh $  that consists of objects satisfying the sheaf property.  Explicitely,  objects of the category
$Sh( \cA^{o} \otimes _k \cB )$  are contravariant  functors $\cA^{o} \otimes _k \cB \to Mod(k)$ that are left exact with respect to both arguments.
The embedding  $Sh \to PSh $ has a 
left adjoint functor (sheafification) 
$$
\tilde \, :  PSh  \to  Sh,
$$
which is exact (\cite{ks}, \S 17.4).  We denote by $\gamma: D(PSh) \to D(Sh) $ the induced functor between the derived categories. The composition
$$\bD( D^{b}_{dg}({\cA})^{o}\otimes_k D^{b}_{dg}({\cB})) \rar{\beta} D(PSh) \rar{\gamma} D(Sh)$$
is {\it not}  fully faithful in general, however, we have the following result.
\begin{lm} \label{key} ({\it cf.} \cite{t},  Th. 8.9)  Let $\bD^+\subset   \bD( D^{b}_{dg}({\cA})^{o}\otimes_k D^{b}_{dg}({\cB})) $ be the full subcategory whose objects are DG modules $M$ such that $\beta(M)$ is bounded from below.
Then the functor
$$S:  \bD^+  \rar{\beta} D^+(PSh) \rar{\gamma} D^+(Sh)$$
is an equivalence of categories.
\end{lm}
\begin{proof} 
 The category $ D^{b}_{dg}({\cA})^{o}\otimes_k D^{b}_{dg}({\cB})$  is the DG quotient of the category $ C^{b}_{dg}({\cA})^{o}\otimes_k C^{b}_{dg}({\cB})$ by the full subcategory whose objects are
 of the form $X^{\cdot} \otimes X^{\prime \cdot}$, where either $X^{ \cdot}$ or $X^{\prime \cdot}$ is acyclic.  It then follows from (\cite{dri}, Theorem 1.6.2) that the functor 
 $$ \beta:  \bD( D^{b}_{dg}({\cA})^{o}\otimes_k D^{b}_{dg}({\cB}))  \to \bD( C^{b}_{dg}({\cA})^{o}\otimes_k C^{b}_{dg}({\cB})) = D(PSh)$$ 
 is fully faithful and that its essential image consists of all  DG-modules $ M \in \bD( C^{b}_{dg}({\cA})^{o}\otimes_k C^{b}_{dg}({\cB}))$ that carry every $X^{\cdot} \otimes X^{\prime \cdot}$ with the above property to an acyclic complex.  
  Identifying the category $ \bD( C^{b}_{dg}({\cA})^{o}\otimes_k C^{b}_{dg}({\cB}))$  with $ D(PSh)$ and observing that the subcategories of acyclic complexes in the homotopy categories  $HoC^{b}_{dg}({\cA})$, $HoC^{b}_{dg}({\cB})$ are generated by short exact sequences
  (\cite{n}, \S 1)
  we exhibit $\bD( D^{b}_{dg}({\cA})^{o}\otimes_k D^{b}_{dg}({\cB})) $ as a full subcategory $\cR \subset D(PSh)$
 whose objects are complexes  $F^{\cdot}$ of presheaves satisfying the following two conditions:
\begin{itemize}
\item For any exact sequence $0\to Z \to Y \to X \to 0$ in $ \cA^{o}$ and any $X' \in  \cB $ the total complex of
\begin{equation}\label{eq2}
 F^{\cdot}(X\otimes X')\to  F^{\cdot}(Y\otimes X') \to F^{\cdot}(Z\otimes X') 
 \end{equation}
is acyclic. 
\item For any  $X \in   \cA^{o}$  and any exact sequence $0\to Z' \to Y' \to X' \to 0$ in $\cB$  the total complex of
$$F^{\cdot}(X\otimes X')\to  F^{\cdot}(X\otimes Y') \to  F^{\cdot}(X\otimes Z') $$
is acyclic.
\end{itemize} 
Observe that, for every $F^{\cdot}\in \cR$ and an exact sequence $0\to Z \to Y \to X \to 0$ in $ \cA^{o}$, we have a long exact sequence of $k$-modules
\begin{equation}\label{eq3}
\cdots H^{m-1}(F^{\cdot}(Z\otimes X')) \to H^{m}(F^{\cdot}(X\otimes X'))\to H^m (F^{\cdot}(Y\otimes X')) \to H^m(F^{\cdot}(Z\otimes X')) \to \cdots 
 \end{equation}
 The equivalence of categories 
$$ \beta:  \bD( D^{b}_{dg}({\cA})^{o}\otimes_k D^{b}_{dg}({\cB}))  \iso \cR \subset D(PSh)$$ 
carries $\bD^+$ to the subcategory  $\cR^+$  of $\cR$ that consists of bounded from below complexes. 

The derived category of sheaves $D(Sh)$ is the quotient of the derived category of presheaves  
by  the subcategory $\cI_{lac}\subset D(PSh)$ of locally (for our Grothendieck topology on  $\cA ^{o} \otimes_k  \cB$)  acyclic complexes  (\cite{bv},  \S 1.11). We shall
prove  that 
\begin{equation}\label{key2}
\cR^+ \subset \cI_{lac}^{\perp},
\end{equation}
 where $\cI_{lac}^{\perp}$ denotes the right orthogonal complement to $\cI_{lac}$ in  $D(PSh)$ (\cite{bv} \S 1.1);  {\it i.e.} 
\begin{equation}\label{keyeq}
Hom_{D(PSh)}(G^{\cdot}, F^{\cdot}) = 0.
\end{equation}
for every $G^{\cdot} \in \cI_{lac}$ and  $F^{\cdot} \in \cR^+$. Without loss of generality we may assume that $F^{\cdot}$ has trivial cohomology in negative degrees: $F^{\cdot}= F^0\to F^1 \to \cdots$. Let $\tilde F^{\cdot} = \tilde F^0\to \tilde F^1 \to \cdots$ be the corresponding complex of sheaves. Since the category of sheaves has enough injective objects (see, e.g. \cite{ks}, Th. 9.6.2, 18.1.6) there exists a complex $I^{\cdot}= I^0 \to I^1 \to \cdots $ of injective sheaves together with a morphism $\tilde F^{\cdot} \to I^{\cdot}$ which is an isomorphism in the derived category of sheaves.  Let us show that the composition $$\delta: F^{\cdot} \to \tilde F^{\cdot} \to  I^{\cdot}$$ is an isomorphism in the derived category of presheaves. 
Indeed, every injective sheaf, viewed as a presheaf, is an object of $\cR$. Thus $I^{\cdot}$ and  $cone(\delta )$ are in $\cR^+$.
Assuming that $cone(\delta )\ne 0$ choose the smallest integer $m$ such that $$0 \ne H^m(cone(\delta )) \in  PSh. $$  Then, there exist an object $X\otimes X'  \in \cA ^{o} \otimes_k  \cB $ and a nonzero element $a\in H^m(cone(\delta ))(X\otimes X')$ . Since the sheafification of $H^m(cone(\delta ))$ is $0$ there exists  a cover $p: Y \otimes Y' \to  X \otimes  X'$ such that
$$0= p^*a \in H^m(cone(\delta ))(Y\otimes Y').$$ Writing $p$ as a composition $$Y \otimes Y' \rar{1\otimes g}   Y\otimes X' \rar{f \otimes 1}    X \otimes X'$$ we may assume $(f\otimes 1)^* a =0$ (otherwise, we replace  $X \otimes X'$ by $ Y \otimes X'$). Let us look at the following fragment of the long 
exact sequence  (\ref{eq3})
applied to $F= cone(\delta)$
and the exact sequence $0\to Z \to Y \rar{f} X\to 0$:   
$$H^{m-1}(cone(\delta ))(Z \otimes X') \to H^{m}(cone(\delta ))(X \otimes X') \to  H^{m}(cone(\delta ))(Y \otimes X'). $$
Since, by our assumption, $ H^{m-1}(cone(\delta )) =0$, it follows that $(f\otimes 1)^*$ is injective and, hence,  $a=0$. This contradiction proves that $cone(\delta )=0$ {\it i.e.},  $\delta$ is a quasi-isomorphism.
 Thus, to complete the proof of (\ref{keyeq}) it suffices to show that 
$$ Hom_{D(PSh)}(G^{\cdot}, I^{\cdot})  =0,$$
for every  $G^{\cdot} \in \cI_{lac}$  and every bounded from below complex of injective sheaves $I^{\cdot}$.  Indeed, every morphism $h:  G^{\cdot} \to I^{\cdot}$ in the derived category is represented by a diagram in $C(PSh( \cA^{o} \otimes _k \cB ))$
$$G^{\cdot} \leftarrow G^{\prime \cdot} \rar{h'} I^{\cdot},$$ 
where the first arrow is a quasi-isomorphism (and, in particular,  $G^{\prime \cdot} \in \cI_{lac}$). If $h'$ is homotopic to $0$ then $h$ is $0$ in the derived category. Thus, it is enough to show that 
 $$Hom_{K(PSh)}(G^{\prime \cdot}, I^{\cdot}) =0, $$
 where $ K(PSh)$ denotes the homotopy category of complexes. We have
 $$ Hom_{K(PSh )}(G^{\prime \cdot}, I^{\cdot}) \iso Hom_{K(Sh)}(\tilde G^{\prime \cdot}, I^{\cdot})
   \iso   Hom_{D(Sh)}(\tilde G^{\prime \cdot}, I^{\cdot}).  $$
 The first arrow is an isomorphism because all terms of the complex $I^{\cdot}$ are sheaves; the second arrow is an isomorphism by (\cite{ks}, Lemma 13.2.4). Finally,  the group $Hom_{D(Sh)}(\tilde G^{\prime \cdot}, I^{\cdot})$ is trivial because the sheafification $\tilde G^{\prime \cdot}$ is $0$ in $D(Sh)$.

To finish the proof of the lemma,  we observe that, for every triangulated category $\cC$ and its full triangulated subcategory
$\cI$, the composition
   $$\cI^{\perp} \to \cC \to \cC/ \cI $$
   is a fully faithful embedding: for every $X,Y \in  \cC$
   $$Hom_{ \cC/ \cI }(X,Y):= \underset{f: X' \to X} {colim} Hom_{\cC} (X', Y),$$
   where the colimit is taken over the filtrant category of pairs $(X' \in \cC, f: X' \to X)$ such that $cone\, f \in \cI$. If $Y\in \cI^{\perp}$, then
   $$ Hom_{\cC}(X,Y) \iso  Hom_{\cC}(X',Y),$$
  and, hence, 
  $$Hom_{ \cC/ \cI }(X,Y)=  Hom_{\cC} (X, Y).$$
   Applying this remark to $\cC= D(PSh)$, $\cI= \cI_{lac}$ and using (\ref{keyeq}) we conclude that the functor $\cR^+ \rar{\gamma}  D(Sh)$ is fully faithful and, hence, so is the composition
   $\bD ^+ \iso \cR^+  \rar{\gamma}  D(Sh)$.  The essential image the functor $\cR^+  \rar{\gamma}  D(Sh)$ coincides with $D^+(Sh)$ because because every complex of injective sheaves viewed as a complex of presheaves is an object of $\cR^+$.    
\end{proof}

\begin{rem}  Applying Lemma \ref{key}  to $k=\bZ$ and $\cA$ being the category of free abelian groups of finite rank  we obtain the following statement:   for every small abelian category $\cB$ 
$$\bD^+(D^{b}_{dg}({\cB}))\iso  D^+(PSh(\cB)) = D^+(Ind(\cB)), $$
where $\bD^+(D^{b}_{dg}({\cB}))$ is the full subcategory of $\bD(D^{b}_{dg}({\cB}))$ that maps to $ D^+(PSh(\cB))$ under the restriction functor (and the ind-completion $Ind(\cB)$ is just another name for $PSh(\cB)$ (\cite{ks}, \S  8.6)). Note the functor
\begin{equation}\label{wrong}
\bD(D^{b}_{dg}({\cB}))\to  D(Ind(\cB)) 
\end{equation}
 is not an equivalence of categories in general.  In fact, the functor  (\ref{wrong})  factors as
 \begin{equation}\label{nice}
  \bD(D^{b}_{dg}({\cB}))\stackrel{\phi}{\iso}  HoC(Ind(\cB))/Ho\overline{C^b_{ac}(\cB) }  \rar{p}  D(Ind(\cB)),
  \end{equation}
where  $Ho\overline{C^b_{ac}(\cB) }$ is the smallest triangulated subcategory of the homotopy category  of acyclic complexes $HoC_{ac}(Ind(\cB))$ that contains {\it finite} acyclic complexes   $HoC^b_{ac}(\cB) $ and closed under arbitrary direct sums; the functor $p$ is the projection
$$ HoC(Ind(\cB))/Ho\overline{C^b_{ac}(\cB) } \to HoC(Ind(\cB))/HoC_{ac}(Ind(\cB)).$$
 The equivalence  $\phi $ can be constructed as follows.  Let $\overline{C^b_{ac}(\cB) }$ be the full subcategory of the DG category  $C(Ind(\cB))$ whose objects are those of    $Ho\overline{C^b_{ac}(\cB) }$.
 The DG quasi-functor $ D^{b}_{dg}({\cB})\to C(Ind(\cB))/\overline{C^b_{ac}(\cB) }$ extends uniquely to
a quasi-functor 
$$\phi_{dg}:  \underrightarrow{ D^{b}_{dg}({\cB})} \to C(Ind(\cB))/\overline{C^b_{ac}(\cB) }$$ 
that commutes with arbitrary direct sums (\cite{bv}, \S1.6.1). Define  $$\phi:= Ho\phi_{dg}.$$ Let us show that $\phi$ is an equivalence of categories. The subcategory   $Ho\overline{C^b_{ac}(\cB) }\subset HoC(Ind(\cB))$ is generated by compact objects (e.g., objects of $HoC^b_{ac}(\cB) $); it follows that the projection  
$HoC(Ind(\cB))\to  HoC(Ind(\cB))/Ho\overline{C^b_{ac}(\cB) }$ carries compact objects of  $HoC(Ind(\cB))$ to compact objects of the quotient category  (\cite{bv}, \S1.4.2). In particular, in the following commutative diagram 
     $$
\def\normalbaselines{\baselineskip20pt
\lineskip3pt  \lineskiplimit3pt}
\def\mapright#1{\smash{
\mathop{\to}\limits^{#1}}}
\def\mapdown#1{\Big\downarrow\rlap
{$\vcenter{\hbox{$\scriptstyle#1$}}$}}
\begin{matrix}
     D^{b}_{dg}({\cB})                  &  =             &     D^{b}_{dg}({\cB}) \cr
   \mapdown{i}         &                      &   \mapdown{j}                           \cr
    \bD(D^{b}_{dg}({\cB}))                   &  \rar{\phi}          &   HoC(Ind(\cB))/Ho\overline{C^b_{ac}(\cB) }    \cr
 \end{matrix}
$$
the image of $j$ consists of compact objects. The same is true for the image of $i$ (\cite{bv}, \S1.7). The functors $i$,$j$ are fully faithful and their images
generate the categories $\bD(D^{b}_{dg}({\cB}))    $ , $HoC(Ind(\cB))/Ho\overline{C^b_{ac}(\cB) } $ respectfully. It follows that $\phi$ is an equivalence of categories.

In general, (e.g., if $\cB$ is the category of finitely generated modules over a finite group )  the projection $p$ is not conservative. However, if the category $\cB$ has {\it finite homological dimension}  the objects of $D^{b}_{dg}({\cB})  $ are compact in $D^{b}_{dg}(Ind({\cB}))  $ \footnote{Indeed, under our finiteness assumption every complex in $D^{b}_{dg}({\cB})$ is quasi-isomorphic to a finite complex of projective objects. Thus it is enough to show that every projective object of $\cB$ is compact in $D(Ind(\cB))$. This is clear because every such object is projective
and compact in $Ind(\cB)$.}  and the above argument proves that (\ref{wrong}) is an equivalence of categories. 
 \end{rem}

 

\begin{cor}\label{cor}
The composition 
\begin{equation}\label{eq40}
S:  {\cT}^+ \rar{\alpha}  \bD( D^{b}_{dg}({\cA})^{o}\otimes_k D^{b}_{dg}({\cB})) \rar{\beta} D(PSh) \rar{\gamma} D(Sh)
\end{equation}
is a fully faithful  embedding. 
\end{cor}

Consider the Yoneda embedding
$$s: Fun(\cA, \cB)\to PSh$$
that takes a functor $F\in Fun(\cA, \cB)$ to the presheaf
$$s(F)(X\times X')= Hom_{\cB}(X', F(X)).$$  
If $F$ is left exact then $s(F)$ is actually a sheaf.

Let $ \cF\in \cT$ be a DG quasi-functor satisfying property $(P)$.  It follows from the definition of  $ \cT^+$ given at the beginning of this section that  $\cF\in \cT^+$. We shall prove that    $S(\cF)\iso s(H^0\cF)$. 
Having in mind applications to Theorem \ref{unbounded} we will actually show a slightly more general statement.  Namely, let us extend the functor (\ref{eq40}) to a larger category:
$$S':  \cT(D^{b}_{dg}({\cA}), D^{+}_{dg}({\cB}))\rar{\alpha'} \bD( D^{b}_{dg}({\cA})^{o}\otimes_k D^{+}_{dg}({\cB})) \rar{\beta '} D(PSh) \rar{\gamma} D(Sh).$$
\begin{lm}\label{key3} Let $ \cF \in \cT(D^{b}_{dg}({\cA}), D^{+}_{dg}({\cB}))$ be a DG quasi-functor
such that $H^i \cF$ is zero for $i<0$ and effaceable for $i>0$. Set $s(F)=s(H^0\cF) \subset Sh\subset D(Sh)$ \footnote{The vanishing of $H^i \cF$ implies that $F$ is left exact and, hence,  $s(F)$ is a sheaf.}.
Then  the complex $S'(\cF) \in D(Sh)$ is canonically quasi-isomorphic to $s(F)$.
\end{lm}
\begin{proof}
By definition, the cohomology presheaves of the complex  $\beta ' \alpha' (\cF) \in D(PSh)$ are given by the formula
$$H^i(\beta' \alpha' \cF)(X\otimes X') = Hom_{D^+(\cB)}(X', Ho(\cF)(X)[i]).$$
Since the negative cohomology of the complex $Ho(\cF)(X)\in D^+(\cB)$ vanishes the same is true for $\beta' \alpha' \cF$ and, thus, we have
$$H^0(\beta' \alpha' \cF)(X\otimes X') = Hom_{D^+(\cB)}(X',  H^0\cF(X)) = s(F).$$
It remains to prove that for every $i>0$ the sheafification of the presheaf $H^i(\beta' \alpha' \cF)$ equals zero. Given an integer $j$ define presheaves $G^{i,j}$ to be
$$ G^{i,j}(X\otimes X') = Hom_{D^+(\cB)}(X', \tau_{\leq j} (Ho(\cF)(X)) [i]).$$
We shall show by induction on $j$ that for every $i>0$ and every $j$ the sheafification of $G^{i,j}$ is $0$.  This would complete the proof since $G^{i,j}$ is isomorphic to $H^i(\beta' \alpha' \cF)(X\otimes X')$ for $j\geq i$. For every $i>0$ and every element $v$
of the group 
$$G^{i,0}(X \otimes X')= Ext^i_{\cB}(X', H^0\cF(X))$$
there exists an epimorphism $Y' \to X'$ such that $v$ is annihilated by the map $$Ext^i_{\cB}(X', H^0\cF(X)) \to Ext^i_{\cB}(Y', H^0\cF(X))$$ (\cite{ks}, Exercise 13.17). This proves that the sheafification of $G^{i,0}$ is $0$.
For the induction step, consider the distinguished triangle $$\tau_{\leq j} (Ho(\cF)(X))\to \tau_{\leq j+1} (Ho(\cF)(X)) \to H^{j+1}\cF(X)[-j-1]$$ and the corresponding long exact sequence
$$ \to  G^{i,j}(X\times X') \to   G^{i,j+1}(X\times X')  \to Hom_{D^b(\cB)}(X',H^{j+1}\cF(X)[-j-1+i]) \to .$$
It follows that $G^{i,j+1}$ fits in a long exact sequence
$$\to G^{i,j} \to G^{i,j+1}\to Ext^{i-j-1}_{\cB}(\cdot, H^{j+1}\cF(\cdot)) \to .$$
 The sheafification of  $G^{i,j}$ is $0$ by the induction assumption, the sheafification of $Ext^{i-j-1}_{\cB}(\cdot, H^{j+1}\cF(\cdot)) $ is $0$ because the functor $H^{j+1}\cF$ is effaceable. Hence, the sheafification of  $G^{i,j+1}$ is $0$ as well.
\end{proof}
Now we are ready to prove the second part of the theorem. Given quasi-functors $\cF, \cG\in \cT$ satisfying property (P) we have by Lemmas \ref{key}, \ref{key3}
\begin{equation}\label{eq13}
Hom_{\cT}(\cF, \cG[i])\iso  Hom_{D(Sh)}(S(\cF), S(\cG)[i]) \iso Ext^i_{Sh}(s(H^0\cF), s(H^0\cG)).
\end{equation}
In particular,  $Hom_{\cT}(\cF, \cG[i])$ is  isomorphic to $Hom_{Fun(\cA, \cB)}(H^0\cF, H^0\cG)$ for $i=0$ (since the functor $s: Fun(\cA, \cB) \to PSh$ is fully faithful) and to $0$ for $i<0$.

To  prove the first part of the theorem we need to recall some facts about DG categories and derived functors.
Let $f: \cC_1 \to \cC_2$  be a  DG functor  between small DG categories. Then the restriction functor 
$f_*: \bD(\cC_2)\to \bD(\cC_1)$
admits a left and a right  adjoint functors  (the derived induction and co-induction functors) 
\begin{equation}\label{ad}
f^*, f^! :  \bD(\cC_1)\to \bD(\cC_2)
\end{equation}
  (\cite{dri}, \S 14.12). In particular, we have the canonical morphisms
\begin{equation}\label{adj}
Id \to f_* f^*, \quad f_* f^! \to Id
\end{equation}
$$Id \to f^! f_*, \quad f^* f_* \to Id.$$
It also follows from the adjunction property that $f^*$ commutes with arbitrary direct sums and that $f^!$ commutes with  arbitrary direct products.
If the the functor $Ho(f): Ho(\cC_1) \to Ho(\cC_2)$ is  fully faithful so is $f_*$ and the first two morphisms in (\ref{adj}) are isomorphisms.   

Recall the definition of the derived DG quasi-functor $RF$ of a left exact functor $F: \cA \to \cB$ from  (\cite{dri}, \S 16). 
Consider the functor 
$$\cT(\cA,   D^{b}_{dg}({\cB}) ) \hookrightarrow \bD( C^{b}_{dg}({\cA})^{o}\otimes_k D^{b}_{dg}({\cB})) \rar{f^*}  \bD( D^{b}_{dg}({\cA})^{o}\otimes_k D^{b}_{dg}({\cB}))$$
induced by the projection
 $$f: C^{b}_{dg}({\cA})^{o}\otimes_k D^{b}_{dg}({\cB}) \to D^{b}_{dg}({\cA})^{o}\otimes_k D^{+}_{dg}({\cB}).$$
 Given a $k$-linear functor $F\in Fun(\cA, \cB ) \to \cT(\cA,   D^{b}_{dg}({\cB}) )$ we define the ``derived functor"
 \begin{equation}\label{derived}
 ``RF" = f^*(F)\in  \bD( D^{b}_{dg}({\cA})^{op}\otimes_k D^{b}_{dg}({\cB})). 
 \end{equation}
The right derived DG quasi-functor $RF:  D^{b}_{dg}({\cA})\to D^{b}_{dg}({\cB})$, if it exists, is an object of $\cT(D^{b}_{dg}({\cA}),  D^{b}_{dg}({\cB}))$ whose image in
$ \bD( D^{b}_{dg}({\cA})^{o}\otimes_k D^{b}_{dg}({\cB}))\supset \cT(D^{b}_{dg}({\cA}),  D^{b}_{dg}({\cB}))$ is
$``RF".$
\begin{lm}\label{key5}
Assume that $F$ is left exact.  Then $``RF" \in \bD^+ \subset \bD( D^{b}_{dg}({\cA})^{op}\otimes_k D^{b}_{dg}({\cB}))$
and the functor 
$ S:  \bD^+ \mono D(Sh)$ takes
$``RF"$ to $s(F)$.
\end{lm}
\begin{proof}
Let $\beta: \bD( D^{b}_{dg}({\cA})^{o}\otimes_k D^{b}_{dg}({\cB})) \to D(PSh)$ be the restriction functor, and let $\gamma:  D(PSh) \to  D(Sh)$ be the sheafification functor. As explained in (\cite{dri}, \S 5) the presheaves $H^i(\beta(``RF"))$ can be computed as follows:
\begin{equation}\label{eq5}
 H^i(\beta(``RF")) (X\otimes X') = colim_Q  \, Hom_{D^b(\cB)}(X',  F(Y^{\cdot})[i]),
 \end{equation}
where the colimit is taken over the filtrant category  $Q$ of pairs $(Y^{\cdot} \in HoC^b_{dg}(\cA), f\in  Hom_{HoC^b_{dg}(\cA)}(X, Y^{\cdot}))$ such that $cone(f)$ is acyclic.  
As the subcategory $Q' \subset Q$ consisting of pairs $(Y^{\cdot}, f)$ with $Y^j=0$  for $j<0$ is cofinal in $Q$, the category $Q$ in the equation (\ref{eq5}) can be replaced by $Q'$.  This proves that $``RF" \in \bD^+$.  Let us show that $\gamma \circ \beta (``RF") \simeq s(F)$. We have $$H^0(\beta(``RF")) (X\otimes X') = colim_{Q'}  \, Hom_{D^b(\cB)}(X',  F(Y^{\cdot}))\simeq $$
$$colim_{Q'}  \, Hom_{D^b(\cB)}(X',  \tau _{\leq 0} F(Y^{\cdot})) \simeq colim_{Q'}  \, Hom_{D^b(\cB)}(X', F(X))=s(F)(X\otimes X').$$
It remains to prove that, for every $i>0$, the sheafification of $H^i(\beta(``RF"))$ is $0$.
Let $s$ be the section of $H^i(\beta(``RF")) (X\otimes X')$ represented by an element $$\tilde s \in Hom_{D^b(\cB)}(X',  F(Y^{\cdot})[i]),$$ where $X\rar{f} Y^0 \to Y^1 \to \cdots$ is an object of $Q'$. Looking at the diagram
    $$
\def\normalbaselines{\baselineskip20pt
\lineskip3pt  \lineskiplimit3pt}
\def\mapright#1{\smash{
\mathop{\to}\limits^{#1}}}
\def\mapdown#1{\Big\downarrow\rlap
{$\vcenter{\hbox{$\scriptstyle#1$}}$}}
\begin{matrix}
  X                             & \rar{f}             &   Y^0                      & \to            & Y^1                  & \to \cdots \cr
   \mapdown{f}         &                      &   \mapdown{Id}        &               & \mapdown{}      &                   \cr
   Y^0                        & \rar{Id}           &  Y^0                       & \to             & 0                     & \to \cdots \cr
 \end{matrix}
$$
we see that the pullback $(f\otimes Id)^* s \in H^i(\beta(``RF")) (Y^0\otimes X')$ is represented by an element of the group $Hom_{D^+(\cB)}(X',  F(Y^0)[i])= Ext^i_{\cB}(X', F(Y^0))$. For any positive $i$ every element of this group is  annihilated by the map
$Ext^i_{\cB}(X', F(Y^0)) \to Ext^i_{\cB}(Y', F(Y^0))$ for some epimorphism $Y' \to X'$.

\end{proof}

Let us prove the first part of the theorem. Let $\cF\in \cT \subset \bD( D^{b}_{dg}({\cA})^{o}\otimes_k D^{b}_{dg}({\cB})) $ be a DG quasi-functor satisfying property (P) together with an isomorphism $F\simeq H^0\cF $.
We need to construct an isomorphism $ \cF \simeq ``RF"$.  By Lemmas \ref{key3}, \ref{key5}  $\cF,  ``RF" $ are objects of $ \bD^+$.  By  Lemma \ref{key}  the functor $S:  \bD^+  \to D(Sh)$ is fully faithful.  Thus, constructing an isomorphism  $ \cF \simeq ``RF"$ is equivalent  to producing
an isomorphism $ S(\cF) \simeq S(``RF")$ in  $ D(Sh)$ which was done in Lemmas \ref{key3}, \ref{key5}.
Theorem \ref{bounded} is proved.

{\bf Proof of theorem \ref{unbounded}.} 
Let 
$\cT^+ \subset \cT:=  \cT(D^{+}_{dg}({\cA}),  D^{+}_{dg}({\cB}))$ be the full triangulated subcategory whose objects are quasi-functors $\cF$ such that, for some integer $n$,  we have
$$ Ho(\cF)( D^{\geq 0}({\cA})) \subset  D^{\geq n}({\cB}). $$  
We shall prove that the composition
$$\cT^+ \mono   \bD( D^{+}_{dg}({\cA})^{o}\otimes_k D^{+}_{dg}({\cB}))  \rar{Res} \bD( D^{b}_{dg}({\cA})^{o}\otimes_k D^{b}_{dg}({\cB})) \to D(Sh)$$
is a fully faithful embedding.  Here $Res$ denotes the restriction functor induced by the embedding
\begin{equation}\label{eq1}
  D^{b}_{dg}({\cA})^{o}\otimes_k D^{b}_{dg}({\cB})  \to D^{+}_{dg}({\cA})^{o}\otimes_k D^{+}_{dg}({\cB}).
 \end{equation}
 
 To show this we need to introduce a bit of notation. 
 If $\cC$ is an abelian category closed under countable direct sums and  
 $$X^0 \rar{\phi _0}  X^1 \rar{\phi_1} X^2 \rar{\phi_2} \cdots  $$  
 is a diagram of complexes  $X^i\in C(\cC)$, we set
$$hocolim\,  X^i =  cone(\bigoplus_i X^i \rar{v}  \bigoplus_i X^i) \in  C(\cC),$$
where $v_{| X^i}: = Id_{X^i}- \phi_i: X^i \to  \oplus_i X^i$.
There is a canonical morphism
$$ hocolim\,  X^i \to colim\, X^i, $$
which is a quasi-isomorphism if countable direct limits in $\cC$ are exact.  If  this is the case, every morphism  $X^{\cdot}\to X'^{\cdot}$ of diagrams that is a term-wise quasi-isomorphism induces a quasi-isomorphism  of the homotopy  colimits \footnote{For the last property, it suffices to assume that  countable direct sums are exact in $\cC$.}.
  Dually, for a category $\cC$ closed under countable products and a diagram   $$\cdots \to X_2 \rar{\phi _1}  X_1 \rar{\phi_0} X_0 ,  $$  we set 
$$holim\,  X_i =  cone(\prod_i X_i \rar{v}  \prod_i X_i)[-1],$$
where $v_i : = p_i - \phi_i p_{i+1}:  \prod X_i  \to  X_i$ and $p_i:  \prod X_i  \to  X_i$ are the projections.

  Let $\bD^f \subset \bD( D^{+}_{dg}({\cA})^{o}\otimes_k D^{+}_{dg}({\cB})) $ be the full subcategory
whose objects  are the covariant DG functors $M:   D^{+}_{dg}({\cA})\otimes_k D^{+}_{dg}({\cB})^{o} \to C(Mod(k))$ such that, for every $X\in  D^{+}_{dg}({\cA})$ and $X' \in D^{+}_{dg}({\cB})$,
the canonical morphism 
 \begin{equation}\label{lim}
 M(X\otimes X') \to holim \, M(X \otimes \tau_{< i} X'),
 \end{equation}
 is a quasi-isomorphism, and,  for every  $X\in  D^{+}_{dg}({\cA})$ and every {\it bounded} $X' \in D^{b}_{dg}({\cB})$,
 the canonical morphism 
 \begin{equation}\label{colim}
  hocolim \, M(\tau_{< i} X \otimes  X') \to  M( X \otimes  X'), 
  \end{equation}
 is a quasi-isomorphism. 
 \begin{rem}\label{rlim} Since countable direct limits are exact in $\cB$, the morphism  $hocolim\,   \tau_{< i} X' \to X'$ is a quasi-isomorphism. 
 Thus,  property (\ref{lim}) is implied by the following: for every integer $n$ and a countable collection $X^{\prime i} \in  D^{\geq n}_{dg}({\cB})$, the morphism 
 $$ M(X\otimes \oplus_i X^{\prime i} ) \to   \prod _i M(X\otimes X^{\prime i} )$$
 is a quasi-isomorphism.
  \end{rem}
  \begin{rem}\label{rcolim}
 Since directed limits are exact in $Mod(k)$  property (\ref{colim})
 is equivalent to the following: for every $X\in  D^{+}_{dg}({\cA})$ and $X' \in {\cB}$, we have
 \begin{equation}\label{colim2}
 colim \, H^0(M(\tau_{<i} X\otimes X')) \iso H^0(M( X\otimes X')).
  \end{equation}
   \end{rem}
 \begin{lm}\label{lm1} 
 The restriction functor
 $$\bD^f \rar{Res} \bD( D^{b}_{dg}({\cA})^{o}\otimes_k D^{b}_{dg}({\cB}))$$
is an equivalence of categories.
\end{lm}
\begin{proof} 
We shall first consider the restriction
$$f_*: \bD( D^{+}_{dg}({\cA})^{o}\otimes_k D^{+}_{dg}({\cB})) \to  \bD( D^{+}_{dg}({\cA})^{o}\otimes_k D^{b}_{dg}({\cB}))$$
and prove that $f^!$ and $f_*$ define mutually inverse equivalences of categories  
\begin{equation}\label{eq21}
\bD( D^{+}_{dg}({\cA})^{o}\otimes_k D^{b}_{dg}({\cB})) \simeq \bD',
\end{equation}
 where $\bD'$ is the full subcategory of  $ \bD( D^{+}_{dg}({\cA})^{o}\otimes_k D^{+}_{dg}({\cB}))$ whose objects are  DG functors $M$ satisfying the property  (\ref{lim}).
   Let us check  that 
  \begin{equation}\label{eq17} 
  f^!(\bD( D^{+}_{dg}({\cA})^{o}\otimes_k D^{b}_{dg}({\cB}))) \subset \bD' .
  \end{equation}
   For every DG functor $f:  \cC_1 \to \cC_2$  between DG categories over a field, the functor $f^!: \bD(\cC_1) \to \bD(\cC_2)$
admits the following concrete description: if $M: \cC_1 \to  C(Mod(k))$ is a  contravariant DG functor and  $X$ is an object of $\cC_2$, we have
\begin{equation}\label{f}
f^!(M)(X)= Hom_{\bD_{dg}(\cC_1)}(f_*^{dg}Hom_{\cC_2}(\cdot, X), M).
\end{equation}
Here $\bD_{dg}(\cC_i)$  denotes the DG derived category of right $\cC_i$-modules, $f_*^{dg}$ the derived restriction functor, and $Hom_{\cC_2}(\cdot, X)$ is the image of $X$  under the Yoneda embedding $\cC_2 \to \bD_{dg}(\cC_2)$.

     We shall prove that 
   $$ hocolim \,  Hom (\cdot \, ,  X \otimes \tau_{< i} X') \to f_*Hom (\cdot \, , X\otimes X') $$ is an isomorphism in $\bD( D^{+}_{dg}({\cA})^{o}\otimes_k D^{b}_{dg}({\cB}))$.
   Together with (\ref{f}) it will imply (\ref{eq17}). 
   By definition of the tensor product of DG categories, for every $Y\otimes Y' \in D^{+}_{dg}({\cA})^{o}\otimes_k D^{b}_{dg}({\cB})$,
   $$Hom(Y\otimes Y',  X\otimes X') =  Hom(Y,  X) \otimes _k Hom(Y', X').$$
   Hence, it is enough to check that the morphism 
    $$ hocolim \,  Hom_{ D^{+}_{dg}({\cB})}(Y', \tau_{< i} Y) \to Hom_{ D^{+}_{dg}({\cB})}(Y' ,Y) $$
    is a quasi-isomorphism, for every $Y' \in D^{b}_{dg}({\cB})$.
     Using the exactness of direct limits in $Mod(k)$  the last assertion is reduced to  
    the formula
    $$colim \, Hom_{D^{b}({\cB})^{o}}(Y',  \tau_{< i} Y) \simeq Hom_{D^{+}({\cB})^{o}}(Y', Y), $$
   which holds because the group $ Hom_{D^{+}({\cB})^{o}}(Y', \tau _{>i} Y)$ is trivial for large $i$.  This proves the assertion (\ref{eq17}). 
   
   Since the functor $Ho(f)$ is fully faithful, we have $$ f_* f^! \iso Id.$$
   Let us check that for every $M\in \bD' $  the canonical morphism $M \to f^! f_* M$ is an isomorphism.
    Set $G=cone(M \to f^! f_* M)$.  As we have just proved
$G$ belongs to $\bD'( D^{b}_{dg}({\cA})^{o}\otimes_k D^{+}_{dg}({\cB}))$. On the other hand,  the isomorphism  $f_* f^! f_* \simeq f_*$  shows that $f_*G$ is $0$. Hence, $G$ is $0$ by  (\ref{lim}).

Next,  consider the DG functor $$ g: D^{b}_{dg}({\cA})^{o}\otimes_k D^{b}_{dg}({\cB}) \to D^{+}_{dg}({\cA})^{o}\otimes_k D^{b}_{dg}({\cB})$$
and show that $g^*$ and $g_*$
 define mutually inverse equivalences  of categories
 \begin{equation}\label{eq20}
   \bD( D^{b}_{dg}({\cA})^{o}\otimes_k D^{b}_{dg}({\cB})) \simeq \bD'',
   \end{equation}
    where $\bD''$ is a full subcategory  of $ \bD( D^{+}_{dg}({\cA})^{o}\otimes_k D^{b}_{dg}({\cB}))$ whose objects are  DG functors $F$ satisfying property  (\ref{colim}).  
 Let us check that 
 \begin{equation}\label{eq18} 
  g^*(\bD( D^{b}_{dg}({\cA})^{o}\otimes_k D^{b}_{dg}({\cB}))) \subset \bD'' .
  \end{equation}
If $M\in \bD( D^{b}_{dg}({\cA})^{o}\otimes_k D^{b}_{dg}({\cB}))$
 is a functor representable by $$Y\otimes Y' \in D^{b}_{dg}({\cA})^{o}\otimes_k D^{b}_{dg}({\cB})$$ then $g^*M$ is represented by the same object  $Y\otimes Y' $ (viewed as an object of $ D^{+}_{dg}({\cA})^{o}\otimes_k D^{b}_{dg}({\cB})$). Hence 
  (\ref{colim2}) is implied by the formula
 $$hocolim \, Hom_{D^{+}_{dg}({\cA})} (Y, \tau _{<i} X) \simeq Hom_{D^{+}_{dg}({\cA})} (Y, X),  \quad Y\in D^{b}_{dg}({\cA}) $$
 proved above (with $\cA$ replaced by $\cB$).
  Since  $g^*$ commutes with arbitrary direct sums and since $\bD( D^{b}_{dg}({\cA})^{o}\otimes_k D^{b}_{dg}({\cB}))$ is the smallest triangulated subcategory that contains representable functors and closed under direct sums,
 $g^*(M)$ is an object of $\bD''( D^{+}_{dg}({\cA})^{o}\otimes_k D^{b}_{dg}({\cB})) $ for every $M$. By  (\ref{colim}) the functor $g_*$ is conservative when restricted to $\bD''$ and the adjoint functor $g^*$ is fully faithful (because $Ho(g)$ is fully faithful).  
 Hence, we have
 $$Id \iso g_* g^*, \quad (g^*g _*)_{ | \bD''} \iso Id.$$
 
 Combining equations (\ref{eq21}) and  $(\ref{eq20})$ we see that the functors $Res$ and $f^! g^*$ define mutually inverse equivalences between the category $\bD^f$ and the category  $\bD( D^{b}_{dg}({\cA})^{o}\otimes_k D^{b}_{dg}({\cB}))$.
 
   \end{proof}
   Consider the composition 
    \begin{equation}\label{eq24} 
   \bD^f \rar{Res} \bD( D^{b}_{dg}({\cA})^{o}\otimes_k D^{b}_{dg}({\cB})) \rar{\beta} D(PSh) \to D(Sh).
    \end{equation}
   Combining Lemmas \ref{key} and \ref{lm1} we get the following.
   
   \begin{cor}\label{keyun} Let $\bD^{f +} \subset   \bD^f $ be the full subcategory whose objects are DG modules $M$ such that $\beta \circ \Res(M)$ is bounded from below.
Then  (\ref{eq24}) induces  an equivalence of categories
$$S:  \bD^{f +} \iso D^+(Sh).$$
\end{cor}

\begin{lm}\label{lm6}  The functor $\cT \mono   \bD( D^{+}_{dg}({\cA})^{o}\otimes_k D^{+}_{dg}({\cB}))$
carries   $\cT^+ $ into $\bD^{f +}$.
\end{lm}
\begin{proof}  Let us show that every $\cF\in \cT$ satisfies property (\ref{rlim}).  By definition of $\cT$, for every $X\in D^{+}_{dg}({\cA})$, there exists $Y\in D^{+}_{dg}({\cB})$ and an isomorphism
$$\cF(X\times ?)\simeq Hom_{ D^{+}_{dg}({\cB})}(? , Y)$$
in the derived category of right $D^{+}_{dg}({\cB})$-modules.  Property (\ref{rlim}) follows because the morphism
$$Hom_{ D^{+}_{dg}({\cB})}(\oplus _i X^{\prime i} , Y)\to \prod_i  Hom_{ D^{+}_{dg}({\cB})}( X^{\prime i} , Y).$$
is a quasi-isomorphism.

 Let us show that every $\cF\in \cT^+$ satisfies the property (\ref{rcolim}). Denote by $Ho(\cF):  D^{+}_{dg}({\cA}) \to D^{+}({\cB})$ the triangulated functor associated with $\cF$. By definition of $Ho(\cF)$ there is a functorial isomorphism
 \begin{equation}\label{eq30} 
  H^0(\cF( X\otimes X')) \simeq Hom _{ D^{+}({\cB})} (X', Ho\cF(X)) 
  \end{equation}
  In order to check (\ref{rcolim})  we will prove a stronger statement:  for every  $X' \in \cB$ 
 the morphism
 \begin{equation}\label{eq11}
 Hom _{ D^{+}({\cB})} (X', Ho\cF(\tau _{<n} X)) \to Hom _{ D^{+}({\cB})} (X', Ho\cF(X))
  \end{equation}
 is an isomorphism for sufficiently large $n$.  By definition of $\cT^+$ we can find an integer $N$ such that
  the functor  $Ho\cF$ carries every object of $D^{> N }({\cA})$ to an object $D^{> 0}(\cB)$.  In particular, for every $n>N$, 
the complex  $Ho\cF(cone(\tau _{<n} X \to X))$ has trivial cohomology in non-positive degrees. Hence, we have
   $$Hom _{ D^{+}({\cB})} (X',  Ho\cF(cone(\tau _{<n} X \to X)))=0.$$

\end{proof}
 
Combining Lemma \ref{lm6} and Corollary \ref{keyun} we get a fully faithful embedding
\begin{equation}\label{eq111}
S: \cT ^+  \mono D(Sh).
\end{equation}
By Lemma \ref{key3}  $S$ carries every quasi-functor $\cF$ satisfying property $(P')$ to $s(H^0\cF)\in Sh$. This proves the second part of Theorem \ref{unbounded}.
For the first part, let  $F\in Fun(\cA, \cB )$ be a $k$-linear functor, and let 
 \begin{equation}\label{derivedun}
 ``RF" \in  \bD( D^{+}_{dg}({\cA})^{o}\otimes_k D^{+}_{dg}({\cB}))
 \end{equation}
 be the ``derived functor" (see (\ref{derived})).
 To complete the proof of Theorem it suffices to show the following.
\begin{lm}
Assume that $F$ is left exact.  Then $``RF" $ is an object of $\bD^{f+} $
and
$ S(``RF")$ is isomorphic to $s(F)$. 
\end{lm}
\begin{proof}
Let us show that $``RF" $ satisfies property (\ref{lim}).  According Remark \ref{rlim} it will suffice to show that,
 for every integer $n$,  $Y^i\in  D^{\geq n}_{dg}({\cB})$ and $ X\in HoC^{+}({\cA})$  
 $$  H^0(``RF" (X\otimes \oplus_i X^{\prime i} ))  \iso   \prod _i   H^0(``RF" (X\otimes X^{\prime i} )). $$
   We have (\cite{dri}, \S 5)
 \begin{equation}\label{eq50}
 H^0(``RF" (X\otimes X')) \simeq  colim_{Q_X}  \, Hom_{D^+(\cB)}(X',  F(Y)),
 \end{equation}
where  $Q_X$ is the filtrant category of pairs $$(Y \in HoC^+_{dg}(\cA), f\in  Hom_{HoC^+_{dg}(\cA)}(X, Y))$$ such that $cone(f)$ is acyclic.  
If $X$ is in $HoC^{\geq n}({\cA})$ the subcategory $Q'_X \subset Q_X$ formed by  pairs $(Y, f)$ with $Y \in  HoC^{\geq n}({\cA}) $   is cofinal in $Q_X$ and, hence,  $Q_X$ in equation (\ref{eq50}) can be replaced by $Q'_X$.  
Thus, it is enough to prove that the category $Q_X$ has the following property: for every countable collection $w_i= (Y_i, f_i)\in Q_X'$, $(i=1,2, \cdots)$, there exists $v\in Q_X$ such that, for every $i$,  the set $Mor_{Q_X}(w_i, v)$ is not empty. 
In fact, the object $$v=( cone(\bigoplus_{i} X \rar{\phi} \bigoplus_{i} Y_i), g),$$ where $\phi_j:  X  \to \bigoplus_{i} Y_i$ equals $f_j-f_{j-1}$ and $g$ is induced by the morphisms $X\rar{f_1} Y_1\mono \bigoplus_{i} Y_i$, does the job. 

Let us show that $``RF" $ satisfies property (\ref{colim}). As we explained in Remark \ref{rcolim} it suffices to show that
$$  colim \, H^0(``RF" (  \tau_{<i} X\otimes X' )) \iso  H^0(``RF" (  X\otimes X' )),$$
for every $X' \in \cB$.
In fact, formula  (\ref{eq50}) with $Q_{ \tau_{\geq i} X}$ replaced by $Q_{ \tau_{\geq i} X}'$ shows that $H^0(``RF" (  \tau_{\geq i} X\otimes X' ))$ is trivial for $i>0$. Hence, the morphism 
$H^0(``RF" (  \tau_{<i} X\otimes X' ))\to H^0(``RF" (  X\otimes X' ))$ is an isomorphism for $i>1$.  This proves that $``RF"$  belongs to $\bD^{f+} $. 

For the second claim,  observe that the restriction  $ Res(``RF" ) \in  \bD( D^{b}_{dg}({\cA})^{o}\otimes_k D^{b}_{dg}({\cB}))$ is the bounded "derived functor"  (\ref{derived}). Thus, we are done by Lemma \ref{key5}.

\end{proof}

{\bf Proof of theorem \ref{Hochschild}.}  Apply  Corollary \ref{cor} and  equation (\ref{eq111}).

\bibliographystyle{alpha}

\begin{thebibliography}{BFM}


\bibitem[BV]{bv} A. Beilinson, V. Vologodsky,
{\it A guide to Voevodsky's motives,}
   Geom. Funct. Anal.  17  (2008),  no. 6, 1709-1787.
   
\bibitem[G]{g}   A. Grothendieck,  {\it Sur quelques points d'alg\`ebre homologique, }  T\^ohoku Math. J. (2)  9  (1957),  119--221. 
   
\bibitem[Dri]{dri} V.Drinfeld, {\it DG quotients of DG categories,}  J. Algebra 272
(2004), no. 2, 643-691.

\bibitem[KS]{ks} M. Kashiwara,  P. Schapira, {\it Categories and sheaves}, A series of comprehensive studies in Mathematics, v. 332, Springer,  (2006).

\bibitem[K1]{k1} B. Keller, {\it On differential graded categories,}  International Congress of Mathematicians. Vol. II,  151--190, Eur. Math. Soc., Z\"urich,  (2006).

\bibitem[K2]{k2}  B. Keller, {\it On the cyclic homology of exact categories,}  J. Pure Appl. Algebra  136  (1999),  no. 1, 1--56.


\bibitem[LV]{lv}  W. Lowen,  M. Van den Bergh,  {\it Hochschild cohomology of abelian categories and ringed spaces,}  Adv. Math.  198  (2005),  no. 1, 172--221.

\bibitem[N]{n} A. Neeman, {\it The derived category of an exact category,} J. Algebra 135
(1990), no. 2, 388-394.

\bibitem[Q]{q} D. Quillen, {\it Higher algebraic K-theory: I}, Higher K-Theories, Lecture Notes in Mathematics, 341 (1972), Springer, pp. 85-147.

\bibitem[T]{t} B. To\"en, {\it The homotopy theory of $dg$-categories and derived Morita theory,}
Invent. Math. 167 (2007), no. 3, 615--667. 

\bibitem[V]{v} V. Vologodsky, {\it The Albanese functor commutes with the Hodge realization,} 
arXiv:0809.2830 (2009). 
\end{thebibliography}

\newcommand\titl[1]{{\newblock{\em #1}}}

\end{document}